\documentclass{article}
\usepackage[OT2,T1]{fontenc}
\usepackage[english,francais]{babel}
\usepackage{lmodern}
\usepackage{amsmath}
\usepackage{amssymb}
\usepackage{mathrsfs}
\usepackage[a4paper,includehead=true,head=1cm,hscale=0.75,top=3cm,bottom=4cm]{geometry}
\newtheorem{theorem}{Theorem}[section]

\newtheorem{e-proposition}[theorem]{Proposition}
\newtheorem{corollary}[theorem]{Corollary}
\newtheorem{e-definition}[theorem]{Definition\rm}

\newtheorem{theoreme}{Th\'eor\`eme}[section]

\newtheorem{proposition}[theoreme]{Proposition}
\newtheorem{corollaire}[theoreme]{Corollaire}

\def\og{\leavevmode\raise.3ex\hbox{$\scriptscriptstyle\langle\!\langle$~}}
\def\fg{\leavevmode\raise.3ex\hbox{~$\!\scriptscriptstyle\,\rangle\!\rangle$}}
\usepackage{amsmath}
\usepackage{mathrsfs}
\newcommand{\Q}{\mathbf{Q}}
\DeclareMathOperator{\Card}{Card}
\newcommand{\fp}[2]{#1_{\mkern-1.5mu#2}}
\newcommand{\ft}[2]{#1_{\mkern-1mu#2}}

\newcommand{\bigensemble}[2]{\left\{ #1 \, ; \, #2 \right\}}
\newcommand{\ensemble}[2]{\{ #1 \, ; \, #2 \}}
\newcommand{\Z}{\mathbf{Z}}
\newcommand{\kbar}{{\mkern1mu\overline{\mkern-1mu{}k\mkern-1mu}\mkern1mu}}
\renewcommand{\P}{\mathbf{P}}
\newcommand{\Br}{{\mathrm{Br}}}
\newcommand{\Gm}{\mathbf{G}_\mathrm{m}}
\renewcommand{\emptyset}{\varnothing}
\newcommand{\uplet}[2]{#1, \mskip2.5mu \ldots \mskip-1mu, \mskip2.5mu #2}
\newcommand{\Srond}{{\mathscr{S}}}
\begin{document}
\selectlanguage{francais}
\title{Principe de Hasse pour les intersections de deux quadriques}
\date{}
\author{\vspace*{.2cm}Olivier Wittenberg\\
\small{}Laboratoire de math\'ematiques, b\^atiment 425,
Universit\'e Paris-Sud, 91405 Orsay, France\\
\small{}\verb!olivier.wittenberg@ens.fr!}
\maketitle

\vspace*{-.8cm}
\selectlanguage{francais}
\begin{center}
{\small 29 novembre 2005}
\end{center}

\vspace*{.2cm}
\begin{abstract}
\selectlanguage{francais}
Admettant l'hypoth\`ese de Schinzel et la finitude des groupes de
Tate-Shafarevich des courbes elliptiques sur les corps de nombres, toute
intersection lisse de deux quadriques dans l'espace projectif de dimension~$n$
satisfait au principe de Hasse si~$n \geqslant 5$.  Le m\^eme r\'esultat vaut
pour $n=4$, c'est-\`a-dire pour les surfaces de del Pezzo de degr\'e~$4$,
lorsque le groupe de Brauer est r\'eduit aux constantes et que la surface est
suffisamment g\'en\'erale.  Les preuves d\'etaill\'ees des r\'esultats
annonc\'es dans cette Note seront publi\'ees ult\'erieurement.
\vskip 0.5\baselineskip

\vspace*{.5cm}
\selectlanguage{english}
\centerline{{\bf Abstract}}
\vskip 0.5\baselineskip
\noindent
{\bf Hasse principle for intersections of two quadrics. }
Assuming Schinzel's hypothesis and the finiteness of Tate-Shafarevich groups
of elliptic curves over number fields, smooth intersections of two quadrics in
$n$\nobreakdash-dimensional projective space satisfy the Hasse principle if~$n
\geq 5$.  The same result holds for $n=4$, \emph{i.e.}, for del Pezzo
surfaces of degree~$4$, provided the Brauer group is reduced to constants and
the surface is sufficiently general.  Detailed proofs of the results announced
herein will be published later on.
\end{abstract}

\bigskip
\selectlanguage{english}
\section*{Abridged English version}

Let~$q_1$ and~$q_2$ be homogeneous quadratic forms in five variables, with coefficients
in a number field~$k$.  Assume that the projective variety~$X \subset \P^4_k$
defined by $q_1=q_2=0$ is a smooth surface; it is
then a del Pezzo surface of degree~$4$, and all del Pezzo surfaces of degree~$4$
are obtained in this way.  The homogeneous polynomial
$f(\lambda,\mu)=\det(\lambda q_1 + \mu q_2) \in k[\lambda,\mu]$
is separable and has degree~$5$.  Let
$\Srond \subset \P^1_k$ be the closed subscheme defined by $f=0$
and let $\{\uplet{t_0}{t_4}\}$ be the set of its $\kbar$-points.
The residue field of a closed point $t \in \Srond$ is denoted~$\kappa(t)$.

For $t \in \Srond$, if $[\lambda:\mu]$ are homogeneous coordinates for~$t$
with $\lambda,\mu\in \kappa(t)$ and if~$L$ is a $\kappa(t)$\nobreakdash-rational
hyperplane of~$\P^4_k$ which does not contain the unique singular point of
the quadric $\lambda q_1 + \mu q_2 = 0$, the discriminant of the restriction
of $\lambda q_1 + \mu q_2$ to the vector space underlying~$L$
does not depend on~$L$.  We shall denote it by
$\varepsilon_t \in \kappa(t)^\star/\kappa(t)^{\star 2}$.

The following proposition gives a practical recipe for computing the cohomological Brauer group
$\Br(X)=H^2_{\text{\'et}}(X,\Gm)$ of~$X$.

\medskip
\begin{e-proposition}
\label{calculbreng}
The group $\Br(X)/\Br(k)$ is a $\Z/2\Z$\nobreakdash-vector space of dimension
$\max(0,n-d-1)$, where $n=\#\bigensemble{t \in
\Srond}{\varepsilon_t \neq 1}$ and~$d$ is the dimension of the
sub-$\Z/2\Z$\nobreakdash-vector space of $k^\star/k^{\star 2}$ spanned
by the norms $N_{\kappa(t)/k}(\varepsilon_t)$ for $t \in \Srond$.
\end{e-proposition}

\medskip
We now state the core technical result of this Note. The reader is referred to~\cite{crelle98}
for more details on Schinzel's hypothesis.

\medskip
\begin{theorem}
\label{thdp4eng}
Assume Schinzel's hypothesis and the finiteness of Tate-Shafarevich groups
of elliptic curves over~$k$.  Assume moreover that $\Br(X)/\Br(k)=0$,
that~$t_0$ is $k$\nobreakdash-rational, that $\varepsilon_t \neq 1$ for each
$t \in \Srond \setminus \{t_0\}$ which has degree at most~$3$ over~$k$,
and that either $\varepsilon_{t_0}=1$ or there exists $t \in \Srond$ such that
the image of~$\varepsilon_{t_0}$ in $\kappa(t)^\star/\kappa(t)^{\star 2}$ is
distinct from~$1$ and from~$\varepsilon_t$.  Then~$X$ satisfies the Hasse principle.
\end{theorem}

\medskip
The most striking consequences of Theorem~\ref{thdp4eng} for the arithmetic
of del Pezzo surfaces of degree~$4$ are summarized below.

\medskip
\begin{corollary}
\label{cordp4eng}
Assume Schinzel's hypothesis and the finiteness of Tate-Shafarevich groups
of elliptic curves over number fields.
The surface~$X$ satisfies the Hasse principle as soon as one of the following
conditions holds:
\selectlanguage{francais}\begin{itemize}\selectlanguage{english}
\item[(i)] the group~$\Gamma$ acts $3$\nobreakdash-transitively on
$\Srond(\kbar)$ (\emph{i.e.} by the symmetric group or by the alternating group);
\item[(ii)] one of the~$t_i$\!'s is $k$\nobreakdash-rational and~$\Gamma$ acts
$2$\nobreakdash-transitively on the other four~$t_i$\!'s;
\item[(iii)] exactly two of the~$t_i$\!'s are $k$\nobreakdash-rational and
$\Br(X)/\Br(k)=0$;
\item[(iv)] all of the~$t_i$\!'s are $k$\nobreakdash-rational and $\Br(X)/\Br(k)=0$.
\end{itemize}
\end{corollary}

\medskip
It was conjectured by Colliot-Th\'el\`ene and Sansuc that~$X$ should satisfy the
Hasse principle as soon as $\Br(X)/\Br(k)=0$ (see~\cite{angers}).
Under Schinzel's hypothesis
and the finiteness of Tate-Shafarevich groups of elliptic curves, Corollary~\ref{cordp4eng}
establishes this conjecture for a very large class of del Pezzo surfaces of degree~$4$.
It should be noted that a sufficiently general del Pezzo surface of degree~$4$
falls into case~(i) above, and that case~(iv) occurs exactly when~$q_1$ and~$q_2$
are simultaneously diagonalizable over~$k$.

\medskip
\begin{corollary}
\label{corint5eng}
Assume Schinzel's hypothesis and the finiteness of Tate-Shafarevich groups of elliptic
curves over number fields.  Smooth intersections of two quadrics in~$\P^n_k$ with $n \geq 5$
satisfy the Hasse principle.
\end{corollary}

\medskip
The conclusion of Corollary~\ref{corint5eng} is a well-known conjecture of Colliot-Th\'el\`ene
and Sansuc (see~\cite[\textsection{}16]{cssint2} and~\cite{harduke}).
Corollary~\ref{corint5eng} is deduced from Corollary~\ref{cordp4eng} by induction on~$n$
and a reduction to hyperplane sections.
When~$n=5$, a monodromy argument is needed to ensure that
the first hypothesis of Corollary~\ref{cordp4eng} is satisfied for a sufficiently
general hyperplane section.

\mathcode`A="7041 \mathcode`B="7042 \mathcode`C="7043 \mathcode`D="7044
\mathcode`E="7045 \mathcode`F="7046 \mathcode`G="7047 \mathcode`H="7048
\mathcode`I="7049 \mathcode`J="704A \mathcode`K="704B \mathcode`L="704C
\mathcode`M="704D \mathcode`N="704E \mathcode`O="704F \mathcode`P="7050
\mathcode`Q="7051 \mathcode`R="7052 \mathcode`S="7053 \mathcode`T="7054
\mathcode`U="7055 \mathcode`V="7056 \mathcode`W="7057 \mathcode`X="7058
\mathcode`Y="7059 \mathcode`Z="705A

\selectlanguage{francais}
\section{Notations}

Soient~$k$ un corps de nombres et~$\kbar$ une cl\^oture alg\'ebrique de~$k$.
Notons~$\Gamma$ le groupe de Galois de~$\kbar$ sur~$k$.  On d\'esigne par
$\P^n_k$ l'espace projectif de dimension~$n$ sur~$k$ et par $(\P^n_k)^\star$
l'espace projectif dual.  Si~$X$ est un sch\'ema, on note~$\kappa(x)$ le corps
r\'esiduel d'un point $x\in X$
et $\Br(X)=H^2_{\text{\'et}}(X,\Gm)$ le groupe de Brauer
cohomologique de~$X$.  Si~$X$ est un $k$\nobreakdash-sch\'ema, on note
abusivement $\Br(X)/\Br(k)$ le conoyau de la fl\`eche naturelle
$\Br(k)\rightarrow \Br(X$); si de plus~$K/k$ est une extension de corps, on
note $X(K)$ l'ensemble des points $K$\nobreakdash-rationnels de~$X$.  Enfin,
on dit qu'une vari\'et\'e~$X$ sur~$k$ \emph{satisfait au principe de Hasse}
lorsque $\prod_{v \in \Omega} X(k_v)\neq\emptyset \Rightarrow
X(k)\neq\emptyset$, o\`u~$\Omega$ d\'esigne l'ensemble des places de~$k$
et~$k_v$ est le compl\'et\'e de~$k$ en~$v$.

\section{Surfaces de del Pezzo de degr\'e~$4$}

Soient~$q_1$ et~$q_2$ deux formes quadratiques homog\`enes en cinq variables,
\`a coefficients dans~$k$.  Notons $X \subset \P^4_k$ la vari\'et\'e
projective d\'efinie par le syst\`eme $q_1=q_2=0$ et supposons que~$X$ soit
une surface lisse.  La surface~$X$ est une surface de del Pezzo de
degr\'e~$4$; r\'eciproquement, toute surface de del Pezzo de degr\'e~$4$ est
isomorphe \`a~$X$ pour des choix appropri\'es de~$q_1$ et~$q_2$.  Le
polyn\^ome $f(\lambda,\mu)=\det(\lambda q_1 + \mu q_2) \in k[\lambda,\mu]$ est
homog\`ene, s\'eparable, de degr\'e~$5$.  Notons $\Srond \subset \P^1_k$ le
sous-sch\'ema ferm\'e d'\'equation $f=0$ et $\{\uplet{t_0}{t_4}\}$ l'ensemble
de ses $\kbar$\nobreakdash-points.

Pour $t \in \Srond$, d\'efinissons une classe $\varepsilon_t \in
\kappa(t)^\star/\kappa(t)^{\star 2}$ comme suit.  Soient~$[\lambda:\mu]$ des
coordonn\'ees homog\`enes pour~$t$, avec $\lambda,\mu\in\kappa(t)$.  Comme~$X$
est lisse, la quadrique d'\'equation $\lambda q_1 + \mu q_2=0$ poss\`ede un
unique point singulier.  Soit~$L$ un hyperplan
$\kappa(t)$\nobreakdash-rationnel de~$\P^4_k$ ne contenant pas ce point.  La
restriction de $\lambda q_1 + \mu q_2$ \`a l'espace vectoriel sous-jacent
\`a~$L$ est une forme quadratique non d\'eg\'en\'er\'ee en quatre variables,
\`a coefficients dans~$\kappa(t)$.  Son discriminant ne d\'epend pas de~$L$;
on le note~$\varepsilon_t$.

\medskip
\begin{proposition}
\label{calculbr}
Le groupe $\Br(X)/\Br(k)$ est un $\Z/2\Z$\nobreakdash-espace vectoriel de
dimension $\max(0,n-d-1)$, o\`u $n=\Card\bigensemble{t \in
\Srond}{\varepsilon_t \neq 1}$ et~$d$ est la dimension du
sous-$\Z/2\Z$\nobreakdash-espace vectoriel de $k^\star/k^{\star 2}$ engendr\'e
par les normes $N_{\kappa(t)/k}(\varepsilon_t)$ pour $t \in \Srond$.
\end{proposition}

\medskip
Le principe de la d\'emonstration est le m\^eme que pour
\cite[Theorem~3.19]{cssint}.

\medskip
Nous aurons \`a consid\'erer la condition suivante:

$(\star)$ le groupe~$\Br(X)/\Br(k)$ est trivial, le point~$t_0$ est
$k$\nobreakdash-rationnel, on a $\varepsilon_t \neq 1$ pour tout $t \in \Srond
\setminus \{t_0\}$ de degr\'e au plus~$3$ sur~$k$, et enfin soit
$\varepsilon_{t_0}=1$, soit il existe $t \in \Srond$ tel que l'image de
$\varepsilon_{t_0}$ dans $\kappa(t)^\star/\kappa(t)^{\star 2}$ soit distincte
de~$1$ et de~$\varepsilon_t$.

\medskip
\begin{theoreme}
\label{thdp4}
Admettons l'hypoth\`ese de Schinzel et la finitude des groupes de
Tate-Shafarevich des courbes elliptiques sur~$k$.  Si la condition~$(\star)$
est v\'erifi\'ee, la surface~$X$ satisfait au principe de Hasse.
\end{theoreme}

\medskip
Nous renvoyons \`a~\cite{crelle98} pour l'\'enonc\'e de l'hypoth\`ese de Schinzel
et pour plus de d\'etails sur le r\^ole qu'elle joue dans ce type de questions.
Le th\'eor\`eme~\ref{thdp4} constitue le c\oe{}ur technique de cette Note.
Avant d'en esquisser la d\'emonstration, citons-en quelques cons\'equences
int\'eressantes pour l'arithm\'etique des surfaces de del Pezzo de
degr\'e~$4$.

\medskip
\begin{corollaire}
\label{cordp4}
Admettons l'hypoth\`ese de Schinzel et la finitude des groupes de
Tate-Shafarevich des courbes elliptiques sur les corps de nombres.  La
surface~$X$ satisfait au principe de Hasse d\`es que l'une des conditions
suivantes est remplie:
\begin{itemize}
\item[(i)] le groupe~$\Gamma$ agit $3$\nobreakdash-transitivement sur
$\Srond(\kbar)$ (\emph{i.e.} par le groupe sym\'etrique ou le groupe
altern\'e);
\item[(ii)] l'un des~$t_i$ est $k$\nobreakdash-rationnel et~$\Gamma$ agit
$2$\nobreakdash-transitivement sur les quatre autres;
\item[(iii)] exactement deux des~$t_i$ sont $k$\nobreakdash-rationnels et
$\Br(X)/\Br(k)=0$;
\item[(iv)] tous les~$t_i$ sont $k$\nobreakdash-rationnels et $\Br(X)/\Br(k)=0$.
\end{itemize}
\end{corollaire}

\medskip
Les surfaces de del Pezzo de degr\'e~$4$ pour lesquelles tous les~$t_i$ sont
$k$\nobreakdash-rationnels sont celles qui sont isomorphes \`a des
intersections de deux quadriques \emph{simultan\'ement diagonales},
c'est-\`a-dire d\'efinies par l'annulation de deux formes quadratiques
diagonales.

Si l'on admet l'hypoth\`ese de Schinzel et la finitude des groupes de
Tate-Shafarevich, le th\'eor\`eme~\ref{thdp4} et le corollaire~\ref{cordp4}
prouvent une grande partie de la conjecture de Colliot-Th\'el\`ene et Sansuc
selon laquelle toute surface de del Pezzo de degr\'e~$4$ telle que
$\Br(X)/\Br(k)=0$ satisfait au principe de Hasse (cf.~\cite{angers}).  Colliot-Th\'el\`ene,
Swinnerton-Dyer et Skorobogatov~\cite[\textsection{}3.2]{css} l'avaient
d\'ej\`a d\'emontr\'ee pour des intersections de deux quadriques
\emph{simultan\'ement diagonales} \og{}suffisamment g\'en\'erales\fg{} (en un sens
explicite), sous l'hypoth\`ese de Schinzel et la finitude des groupes de
Tate-Shafarevich.  Le corollaire~\ref{cordp4}, restreint au cas~(iv),
g\'en\'eralise ce r\'esultat.  Les autres cas sont enti\`erement nouveaux.  Il
est \`a noter qu'une surface de del Pezzo de degr\'e~$4$ \og{}suffisamment
g\'en\'erale\fg{} v\'erifie~(i).

La preuve du corollaire~\ref{cordp4} \`a partir du th\'eor\`eme~\ref{thdp4}
repose sur les trois faits suivants.  Tout d'abord, il r\'esulte de la
d\'emonstration de la proposition~\ref{calculbr} que $\prod_{t \in
\Srond}N_{\kappa(t)/k}(\varepsilon_t)=1$ dans $k^\star/k^{\star 2}$.  Ensuite,
Coray~\cite{coray} a \'etabli que s'il existe une extension finie $K/k$ de
degr\'e impair telle que $X(K)\neq\emptyset$, alors $X(k)\neq\emptyset$.  En
particulier, le cas~(i) se ram\`ene imm\'ediatement au cas~(ii).  Enfin, s'il
existe $t \in \Srond$ de degr\'e impair sur~$k$ tel que $\varepsilon_t=1$ et
si $\Br(X)/\Br(k)=0$, alors~$X$ satisfait au principe de Hasse.  En effet, vu
le r\'esultat de Coray, il suffit de montrer que $X \otimes_k \kappa(t)$
satisfait au principe de Hasse sur~$\kappa(t)$; or les hypoth\`eses
entra\^inent d'une part que l'obstruction de Brauer-Manin \`a l'existence d'un
point $K$\nobreakdash-rationnel s'\'evanouit et d'autre part que $X \otimes_k
\kappa(t)$ contient un pinceau de coniques, et Salberger~\cite{salb} a
d\'emontr\'e que les surfaces de del Pezzo de degr\'e~$4$ admettant un pinceau
de coniques satisfont au principe de Hasse d\`es que l'obstruction de
Brauer-Manin ne s'y oppose pas.

\medskip
\noindent\textbf{Principe de la d\'emonstration du th\'eor\`eme~\ref{thdp4}. }
Le point de d\'epart est la construction sugg\'er\'ee par Swinnerton-Dyer
dans~\cite[\textsection{}6]{bsd} et que l'on peut r\'esumer comme suit.  Pour
$i \in \{\uplet{0}{4}\}$, notons $\fp{P}i \in \P^4(\kbar)$ l'unique point
singulier de la quadrique d'\'equation $\lambda q_1 + \mu q_2=0$, o\`u
$[\lambda:\mu]$ sont des coordonn\'ees homog\`enes pour $t_i \in \P^1(\kbar)$.
D'apr\`es $(\star)$, le point~$\fp{P}0$ est $k$-rationnel.  L'hyperplan
\mbox{$H \subset \P^4_k$} contenant $\uplet{\fp{P}1}{\fp{P}4}$ est donc lui
aussi $k$-rationnel.  Posons $H^0=H \setminus ((X \cap H) \cup
\{\uplet{\fp{P}1}{\fp{P}4}\})$ et \mbox{$Y=\ensemble{(x,h) \in X \times_k
H^0}{x \in \ft{T}h Q_h}$}, o\`u~$Q_h$ d\'esigne l'unique quadrique
contenant~$X$ et~$h$ et~$\ft{T}h Q_h$ est l'espace tangent \`a~$Q_h$ en~$h$.
La seconde projection $\pi \colon Y \rightarrow H^0$ est propre et plate; on
peut donc parler de fibration en courbes de genre~$1$.  Cette fibration jouit
des propri\'et\'es remarquables suivantes: la jacobienne de sa fibre
g\'en\'erique admet un point rationnel d'ordre~$2$, sa base est une
vari\'et\'e $k$\nobreakdash-rationnelle et ses fibres sont des sections
hyperplanes de~$X$.  Par ailleurs, la premi\`ere projection
$Y \rightarrow X$ admet une section rationnelle.

Swinnerton-Dyer et ses collaborateurs ont d\'evelopp\'e une technique pour
\'etudier les points rationnels de certaines surfaces fibr\'ees en courbes de
genre~$1$ au-dessus de~$\P^1_k$, lorsque la jacobienne de la fibre
g\'en\'erique admet un point rationnel d'ordre~$2$ (cf.~\cite{egloff},
\cite{css}, \cite{bsd}, \cite{ctbsd}).  \`A quelques hypoth\`eses techniques
pr\`es, la variante \'etudi\'ee
dans~\cite[\textsection{}1--\textsection{}5]{bsd} et~\cite{ctbsd} s'applique
\`a la restriction de~$\pi$ au-dessus d'une droite suffisamment g\'en\'erale
de~$H$.  Reprenant les arguments de~\cite{bsd} et de~\cite{ctbsd} en termes de mod\`eles
de N\'eron, nous
obtenons d'abord un \'enonc\'e dans lequel ces hypoth\`eses techniques n'apparaissent
plus.  Il en r\'esulte, dans la situation qui nous int\'eresse, que pour
\'etablir l'existence d'un point rationnel sur~$Y$ (et donc sur~$X$) en
admettant l'hypoth\`ese de Schinzel et la finitude des groupes de
Tate-Shafarevich des courbes elliptiques sur~$k$, il suffit d'exhiber une
droite suffisamment g\'en\'erale $L \subset H$ et un point $h \in L(k)$ tels
que la surface $\pi^{-1}(L)$ remplisse une condition du type
\og{}condition~$(D)$\fg{} (cf. \cite[\textsection{}4]{css}) et que la fibre
$\pi^{-1}(h)$ admette un $k_v$\nobreakdash-point pour tout $v \in \Omega$.

Pour les vari\'et\'es munies d'un morphisme vers~$\P^1_k$ dont chaque fibre
singuli\`ere poss\`ede une composante irr\'eductible de multiplicit\'e~$1$
d\'eploy\'ee par une extension ab\'elienne du corps de base, il est connu
qu'en l'absence d'obstruction de Brauer-Manin \`a l'existence d'un point
rationnel, on peut trouver des fibres au-dessus de points rationnels de~$\P^1_k$
comportant un $k_v$\nobreakdash-point pour
tout $v \in \Omega$, si l'on admet
l'hypoth\`ese de Schinzel (cf.~\cite[Theorem~1.1]{crelle98}).  Nous
g\'en\'eralisons ce th\'eor\`eme au cas de vari\'et\'es fibr\'ees au-dessus
de~$\P^n_k$; la d\'emonstration utilise de mani\`ere essentielle un r\'esultat
r\'ecent d'Harari (cf.~\cite{harari}).  Dans la situation consid\'er\'ee ici,
on obtient ainsi l'existence de $h_0 \in H^0(k)$ tel que la fibre
$\pi^{-1}(h_0)$ soit lisse et admette un $k_v$\nobreakdash-point pour tout $v
\in \Omega$.

Nous d\'efinissons ensuite une version de la
\og{}condition~$(D)$\fg{} pour la fibration~$\pi$, dite \emph{condition~$(D)$ g\'en\'erique},
adapt\'ee \`a l'ensemble des points de~$H^0$ de codimension~$\leqslant 1$.
\`A l'aide du th\'eor\`eme d'irr\'eductibilit\'e de Hilbert, nous prouvons
qu'elle implique l'existence de droites $L \subset H$ passant par $h_0$ et
telles qu'une version arithm\'etique de la \og{}condition~$(D)$\fg{} soit
satisfaite pour la surface $\pi^{-1}(L)$.  Une \'etude d\'elicate de
l'arithm\'etique des fibres singuli\`eres de~$\pi$ permet enfin d'\'etablir
que sous l'hypoth\`ese de $k$\nobreakdash-rationalit\'e de~$t_0$
(qui est de toute mani\`ere n\'ecessaire \`a la d\'efinition de~$\pi$),
la condition~$(D)$ g\'en\'erique \'equivaut \`a~$(\star)$.

\section{Intersections lisses de deux quadriques dans~$\P^n_k$ pour $n \geqslant 5$}

\medskip
\begin{theoreme}
\label{intp5}
Admettons l'hypoth\`ese de Schinzel et la finitude des groupes de
Tate-Shafarevich des courbes elliptiques sur les corps de nombres.  Soit $n
\geqslant 5$.  Toute intersection lisse de deux quadriques dans~$\P^n_k$
satisfait au principe de Hasse.
\end{theoreme}

\medskip
La conclusion de ce th\'eor\`eme est une conjecture de Colliot-Th\'el\`ene et
Sansuc (cf.~\cite[\textsection{}16]{cssint2}, \cite{harduke}).
Elle fut notamment \'etablie par Mordell (1959)
pour $n \geqslant 12$ et $k=\Q$, par Swinnerton-Dyer (1964) pour $n \geqslant
10$ et $k=\Q$ et par Colliot-Th\'el\`ene, Sansuc et Swinnerton-Dyer (1987)
pour $n \geqslant 8$ (cf.~\cite{cssint}).

Nous d\'eduisons le th\'eor\`eme~\ref{intp5} du corollaire~\ref{cordp4},
cas~(i), par la m\'ethode des fibrations (cf.~\cite{cssint}, \cite{skofib}).  Le r\'esultat
suivant, appliqu\'e \`a $n=5$, constitue le point cl\'e de la d\'emonstration;
il montre que l'action de~$\Gamma$ sur l'ensemble~$\Srond(\kbar)$ associ\'e
\`a une section hyperplane suffisamment g\'en\'erale d'une intersection lisse
et de dimension~$3$ de deux quadriques dans~$\P^5_k$ est
$3$\nobreakdash-transitive (et m\^eme $5$\nobreakdash-transitive).

\medskip
\begin{proposition}
Soit $X \subset \P^n_k$ une intersection lisse de deux quadriques, de
codimension~$2$ dans~$\P^n_k$.  Notons $D\simeq \P^1_k$ le pinceau des
quadriques contenant~$X$ et $Z=\ensemble{(q,h) \in D \times_k
(\P^n_k)^\star}{q \cap h \text{ est singulier}}$.  La seconde projection $Z
\rightarrow (\P^n_k)^\star$ fait de~$Z$ un rev\^etement int\`egre, fini et
plat de~$(\P^n_k)^\star$.  Son groupe de monodromie est isomorphe au groupe
sym\'etrique $\mathfrak{S}_n$.
\end{proposition}

\section*{Remerciements}

Mes plus vifs remerciements vont \`a Jean-Louis Colliot-Th\'el\`ene, qui me
sugg\'era le probl\`eme et me fit g\'en\'ereusement part de ses id\'ees pour
le r\'esoudre.  Les r\'esultats annonc\'es ici doivent beaucoup \`a Sir Peter
Swinnerton-Dyer, puisqu'il est aussi bien \`a l'origine de la construction
utilis\'ee dans la preuve du th\'eor\`eme~\ref{thdp4} que de la technique
employ\'ee pour \'etudier l'arithm\'etique des pinceaux de courbes de
genre~$1$ (cf.~\cite{egloff}).  Je tiens \'egalement \`a le remercier pour ses
encouragements et pour son aide \`a la compr\'ehension d'un d\'etail important
de~\cite[\textsection{}6]{bsd}.  Enfin, je remercie David Harari d'avoir
r\'edig\'e et de m'avoir transmis~\cite{harari}, dont le th\'eor\`eme
principal joue un r\^ole capital dans la d\'emonstration du
th\'eor\`eme~\ref{thdp4}.

\end{document}